\documentclass[11pt]{article}
\usepackage{a4wide}

\newcommand{\ba}{\begin{array}}
\newcommand{\ea}{\end{array}}
\newcommand{\be}{\begin{equation}}
\newcommand{\ee}{\end{equation}}
\newcommand{\la}{\label}
\newcommand{\bea}{\begin{eqnarray}}
\newcommand{\eea}{\end{eqnarray}}
\newcommand{\ch}{\choose}
\renewcommand{\a}{\alpha}
\renewcommand{\b}{\beta}
\renewcommand{\c}{\gamma}
\newcommand{\G}{\Gamma}
\renewcommand{\L}{L_n^{(\a)}(x)}
\newcommand{\KL}{L_n^{\a,M}(x)}
\newcommand{\gL}{L_n^{\a,M,N}(x)}
\renewcommand{\P}{P_n^{(\a,\b)}(x)}
\newcommand{\GP}{P_n^{\a,\b,M,N}(x)}
\newcommand{\SP}{P_n^{(\a,\a)}(x)}
\newcommand{\SGP}{P_n^{\a,\a,M,M}(x)}
\renewcommand{\l}{\left}
\renewcommand{\r}{\right}
\newcommand{\set}[1]{\left\{#1\right\}_{n=0}^{\infty}}
\newcommand{\hyp}[5]{\mbox{}_{#1}F_{#2}
\left(\left.\begin{array}{c}#3\\#4\end{array}\right|#5\right)}
\newcommand{\ndots}{n=0,1,2,\ldots}
\newcommand{\n}{\nonumber}
\newcommand{\nn}{\nonumber \\}
\newcommand{\ds}{\displaystyle}
\newcounter{stelling}
\newcommand{\st}[1]{\par\vspace{0.5cm}\refstepcounter{stelling}
{\bf Theorem \thestelling.} {\sl #1}\par\vspace{0.5cm}}

\begin{document}

\begin{center}
{\Large THE SEARCH FOR\\ DIFFERENTIAL EQUATIONS\\ FOR\\ ORTHOGONAL POLYNOMIALS\\
BY USING COMPUTERS}

\vspace{1cm}

{\large Roelof Koekoek}

\vspace{1cm}

Delft University of Technology\\ Faculty of Technical Mathematics and Informatics\\
Mekelweg 4\\ 2628 CD Delft\\ The Netherlands\\ E-mail~:
koekoek@twi.tudelft.nl
\end{center}

\vspace{2cm}

\begin{abstract}
We look for differential equations of the form
$$\sum_{i=0}^{\infty}c_i(x)y^{(i)}(x)=\lambda_ny(x),$$
where the coefficients $\l\{c_i(x)\r\}_{i=0}^{\infty}$ do not depend on $n$,
for the generalized Jacobi polynomials $\set{\GP}$ found by T.H.~Koornwinder
in 1984 and for generalized Laguerre polynomials $\set{\gL}$ which are
orthogonal with respect to an inner product of Sobolev type.

We introduce a method which makes use of computeralgebra packages like Maple
and Mathematica and we will give some preliminary results.
\end{abstract}

\newpage

\section{Introduction}

In 1984 T.H. Koornwinder (see \cite{Koorn}) introduced the polynomials $\set{\GP}$ which are
orthogonal on the interval $[-1,1]$ with respect to the weight function
$$\frac{\G(\a+\b+2)}{2^{\a+\b+1}\G(\a+1)\G(\b+1)}(1-x)^{\a}(1+x)^{\b}
+M\delta(x+1)+N\delta(x-1).$$
These generalized Jacobi polynomials generalize the Legendre type and the
Jacobi type polynomials found by H.L.~Krall. See for instance \cite{Krall},
\cite{Krall2} and \cite{Krall3}.

As a limit case Koornwinder also found generalized Laguerre polynomials
$\set{\KL}$ which are orthogonal on the interval $[0,\infty)$ with respect
to the weight function
$$\frac{1}{\G(\a+1)}x^{\a}e^{-x}+M\delta(x).$$

It is well-known that these generalized Jacobi and generalized Laguerre polynomials
satisfy a linear second order differential equation with coefficients
depending on $n$, but of bounded degree. See for instance \cite{Koorn}.

Now we are looking for differential equations of the form
$$\sum_{i=0}^{\infty}c_i(x)y^{(i)}(x)=\lambda_ny(x),$$
where the coefficients $\l\{c_i(x)\r\}_{i=0}^{\infty}$ are independent of the
degree $n$. In \cite{Krall}, \cite{Krall2} and \cite{Krall3} several special
cases were treated. Later more special cases were found by A.M.~Krall and
L.L.~Little\-john. See \cite{ClassI}, \cite{ClassII} and \cite{Lit_Krall}.

These differential equations we are looking for are studied in spectral theory.
A survey of orthogonal polynomials and spectral theory was given in \cite{Conj}.
A large number of references can be found there.

In \cite{DV} J.~Koekoek and R.~Koekoek found a differential equation of the
above type for the generalized Laguerre polynomials $\set{\KL}$ for all $\a>-1$.
This result generalizes the results found for special cases of the Laguerre
type polynomials. In fact we showed that this differential equation is of
infinite order in general, but reduces to finite order for nonnegative integer
values of $\a$. See \cite{DV} for more details and section 3 of this paper
for more results concerning the coefficients of this differential equation.

In \cite{SIAM} R.~Koekoek and H.G.~Meijer introduced the polynomials $\set{\gL}$
which are orthogonal with respect to the following (Sobolev) inner product~:
$$<f,g>\;=\frac{1}{\G(\a+1)}\int\limits_0^{\infty}x^{\a}e^{-x}f(x)g(x)dx+
Mf(0)g(0)+Nf'(0)g'(0),$$
where $\a>-1$, $M\ge 0$ and $N\ge 0$. These polynomials $\set{\gL}$ generalize
the polynomials $\set{\KL}$ since $L_n^{\a,M,0}(x)=\KL$.

In this paper we deal with the problem of finding a differential equation of the
above type for both the generalized Jacobi polynomials $\set{\GP}$ and the
polynomials $\set{\gL}$.

The differential equation found in \cite{DV} was computed by hand, without the
help of computers. This is nearly impossible for these other cases. We need
computers to handle the very huge expressions we have to deal with.

\section{The method}

In this section we will describe the method we use to discover a linear
differential equation satisfied by polynomials which are defined as a linear
combination of classical orthogonal polynomials and their derivatives. These
polynomials depend on $x$, the parameters of the classical orthogonal polynomials
and on some extra parameters, say $M$ and $N$.

As an example we look at the polynomials $\set{\gL}$ which are defined by
$$\gL=A_0\L+A_1\frac{d}{dx}\L+A_2\frac{d^2}{dx^2}\L,$$
where the coefficients $A_0$, $A_1$ and $A_2$ depend on $n$, $\a$, $M$ and $N$.
Moreover, $A_0$, $A_1$ and $A_2$ are linear combinations of $1$, $M$, $N$ and
$MN$. These generalized Laguerre polynomials are defined in such a way that
$L_n^{\a,0,0}(x)=\L$. See \cite{SIAM} and section 4 of this paper for more details.

Now we start from the well-known (second order) differential equation for the
classical orthogonal polynomials, which is
$$xy''(x)+(\a+1-x)y'(x)+ny(x)=0$$
in the Laguerre case. If we substitute the new generalized orthogonal
polynomials into the left-hand side of this classical differential
equation we get an expression which can be seen as a polynomial in $M$ and $N$.
Now we add a linear combination of $M$, $N$ and $MN$ depending on $y$, $y'$,
$y''$, etc. to the left-hand side of this classical differential equation
and set this expression equal to zero. We hope that this is the new differential
equation we are looking for.

Our example leads to the differential equation
\bea & &M\sum_{i=0}^{\infty}a_i(x)y^{(i)}(x)+
N\sum_{i=0}^{\infty}b_i(x)y^{(i)}(x)+{}\nn
& &\hspace{1cm}{}+MN\sum_{i=0}^{\infty}c_i(x)y^{(i)}(x)+
xy''(x)+(\a+1-x)y'(x)+ny(x)=0.\n\eea

Now we substitute the new orthogonal polynomials into this new differential
equation and then the left-hand side becomes a polynomial in $M$ and $N$
with coefficients involving the classical orthogonal polynomials.
Since this differential equation must be valid for all possible values of $M$
and $N$, and $M$ and $N$ are supposed to be independent, all coefficients of
this polynomial must be equal to zero. This gives us several equations in
terms of unknown coefficients and of the well-known classical orthogonal
polynomials. We remark that in case of dependence of $M$ and $N$ we have to
deal with a completely different problem, since then we have a polynomial
in one variable instead of two. We want the coefficients of $y'$, $y''$, etc. to be independent
of the degree $n$. So, if we substitute small values of $n$ into these
equations we get rather simple equations for the coefficients we are looking
for. In general, these equations are huge expressions and are very difficult
to handle. But here we use computers to do these huge calculations. The
mathematics used is very simple, but the formulas are too big to do this by
hand. The calculations we did here were done by using Mathematica and sometimes
Maple. Both computeralgebra packages can do these calculations easily, but
they use a lot of computermemory. On a workstation with 8 megabytes of
internal memory some calculations took several hours.

\section{The infinite order Laguerre differential equation}

In \cite{DV} the following theorem was proved without the use of computers~:

\st{For $M>0$ the polynomials $\set{\KL}$ satisfy a unique differential equation
of the form
\be\la{dv}M\sum_{i=0}^{\infty}a_i(x)y^{(i)}(x)+xy''(x)+(\a+1-x)y'(x)+ny(x)=0,\ee
where $\l\{a_i(x)\r\}_{i=0}^{\infty}$ are continuous functions on the real
line and $\l\{a_i(x)\r\}_{i=1}^{\infty}$ are independent of $n$.

Moreover, the functions $\l\{a_i(x)\r\}_{i=0}^{\infty}$ are polynomials given by
\be\la{coeff}\l\{\ba{l}\ds a_0(x)={n+\a+1 \ch n-1}\\ \\
\ds a_i(x)=\frac{1}{i!}\sum_{j=1}^i(-1)^{i+j+1}{\a+1 \ch j-1}{\a+2 \ch i-j}(\a+3)_{i-j}x^j,
\;i=1,2,3,\ldots.\ea\r.\ee}

For $\a\ne 0,1,2,\ldots$ we have degree$[a_i(x)]=i,\;i=1,2,3,\ldots$. This implies
that if $M>0$ the differential equation (\ref{dv}) is of infinite order in that
case. For nonnegative integer values of $\a$ we have
$$\l\{\ba{ll}\mbox{degree}[a_i(x)]=i, & i=1,2,3,\ldots,\a+2\\ \\
\mbox{degree}[a_i(x)]=\a+2, & i=\a+3,\a+4,\a+5,\ldots,2\a+4\\ \\
a_i(x)=0, & i=2\a+5,2\a+6,2\a+7,\ldots.\ea\r.$$
This implies that for nonnegative integer values of $\a$ and $M>0$ the differential
equation (\ref{dv}) is of order $2\a+4$.

This differential equation was found by setting $y(x)=\KL$ and by substituting
small values of $n$ in (\ref{dv}). Since the coefficients $\l\{a_i(x)\r\}_{i=1}^{\infty}$
are independent of $n$ this gives us $a_1(x), a_2(x), a_3(x), \ldots$ explicitly.
Then the general form of $a_i(x)$ was guessed and the result was proved.
Note that $a_0(x)$ does not depend on $x$, but does depend on $n$. From the
three special cases ($\a=0$, $\a=1$ and $\a=2$) found by A.M.~Krall and
L.L.~Littlejohn the general form of this coefficient could be guessed rather
easily and be proved too.

Later we discovered that the coefficients $\l\{a_i(x)\r\}_{i=1}^{\infty}$
have the following interesting property.

\st{The coefficients $\l\{a_i(x)\r\}_{i=1}^{\infty}$ of the differential
equation given by (\ref{dv}) and (\ref{coeff}) satisfy
\be\la{suma}\sum_{i=1}^{\infty}a_i(x)=-\frac{\sin\pi\a}{\pi}\frac{x}{(\a+2)(\a+3)}
\hyp{1}{1}{1}{\a+4}{-x},\;\a>-1.\ee

For nonnegative integer values of $\a$ we have~:
\be\la{suma2}\sum_{i=k}^{\infty}{i \ch k}a_i(x)=(-1)^{\a+k}a_k(-x),\;k=1,2,3,\ldots.\ee}

Note that this theorem implies for nonnegative integer values of $\a$~:
$$\sum_{i=1}^{\infty}a_i(x)=0\;\mbox{ and }\;\sum_{i=1}^{\infty}ia_i(x)=(-1)^{\a+1}x.$$

{\bf Proof.} First we prove (\ref{suma}). Changing the order
of summation we find
\bea\sum\limits_{i=1}^{\infty}a_i(x)&=&\sum\limits_{i=1}^{\infty}\sum\limits_{j=1}^i
\frac{(-1)^{i+j+1}}{i!}{\a+1 \ch j-1}{\a+2 \ch i-j}(\a+3)_{i-j}x^j\nn
&=&\sum\limits_{j=1}^{\infty}(-1)^{j+1}{\a+1 \ch j-1}x^j\sum\limits_{i=j}
^{\infty}\frac{(-1)^i}{i!}{\a+2 \ch i-j}(\a+3)_{i-j}\nn
&=&\sum\limits_{j=1}^{\infty}(-1)^{j+1}{\a+1 \ch j-1}x^j\sum\limits_{i=0}
^{\infty}\frac{(-1)^{i+j}}{(i+j)!}{\a+2 \ch i}(\a+3)_i\nn
&=&-\sum\limits_{j=1}^{\infty}{\a+1 \ch j-1}x^j\sum\limits_{i=0}^{\infty}
\frac{(-1)^i}{(i+j)!}{\a+2 \ch i}(\a+3)_i.\n\eea
Now we use the well-known summation formula
\be\la{sum}\hyp{2}{1}{a,b}{c}{1}=\frac{\G(c-a-b)
\G(c)}{\G(c-a)\G(c-b)},\;c-a-b>0,\;c\ne 0,-1,-2,\ldots\ee
to find
\bea\sum_{i=0}^{\infty}\frac{(-1)^i}{(i+j)!}{\a+2 \ch i}(\a+3)_i&=&
\frac{1}{j!}\sum_{i=0}^{\infty}\frac{(-\a-2)_i(\a+3)_i}{(j+1)_ii!}\nn
&=&\frac{1}{j!}\hyp{2}{1}{-\a-2,\a+3}{j+1}{1}\nn
&=&\frac{\G(j)}{\G(j+\a+3)\G(j-\a-2)},
\;j=1,2,3,\ldots.\n\eea
Hence
\bea\sum_{i=1}^{\infty}a_i(x)&=&-\sum_{j=1}^{\infty}{\a+1 \ch j-1}
\frac{\G(j)}{\G(j+\a+3)\G(j-\a-2)}x^j\nn
&=&-x\sum_{j=0}^{\infty}{\a+1 \ch j}\frac{\G(j+1)}
{\G(j+\a+4)\G(j-\a-1)}x^j\nn
&=&-x\sum_{j=0}^{\infty}(-1)^j\frac{(-\a-1)_j}{\G(\a+4)(\a+4)_j
\G(-\a-1)(-\a-1)_j}x^j\nn
&=&-\frac{x}{\G(\a+4)\G(-\a-1)}\sum_{j=0}^{\infty}
\frac{(-x)^j}{(\a+4)_j}\nn
&=&-\frac{x}{\G(\a+4)\G(-\a-1)}\hyp{1}{1}{1}{\a+4}{-x}.\n\eea
Finally we use
$$\frac{1}{\G(z)\G(1-z)}=\frac{\sin\pi z}{\pi}$$
to obtain
\bea\frac{1}{\G(\a+4)\G(-\a-1)}&=&\frac{1}{(\a+2)(\a+3)}
\frac{1}{\G(\a+2)\G(-\a-1)}\nn
&=&\frac{1}{(\a+2)(\a+3)}\frac{\sin\pi(\a+2)}{\pi}=
\frac{1}{(\a+2)(\a+3)}\frac{\sin\pi\a}{\pi}.\n\eea
This proves (\ref{suma}).

To prove (\ref{suma2}) we take $\a\in\{0,1,2,\ldots\}$ and start with
$$\sum_{i=k}^{\infty}\frac{i!}{(i-k)!}a_i(x)=\sum_{i=k}^{\infty}\sum_{j=1}^i
\frac{(-1)^{i+j+1}}{(i-k)!}{\a+1 \ch j-1}{\a+2 \ch i-j}
(\a+3)_{i-j}x^j.$$
Now we use
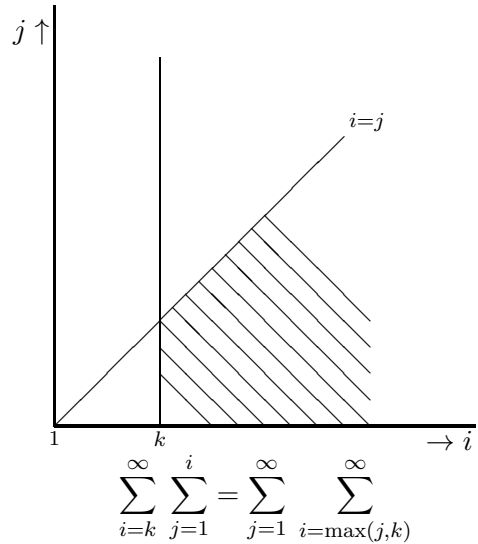
\renewcommand{\ss}{\scriptstyle}
\setlength{\unitlength}{0.7cm}
\begin{center}
\begin{picture}(8,8)\thicklines
\put(0,0){\line(1,0){8}} \put(0,0){\line(0,1){8}}
\thinlines
\put(0,0){\line(1,1){5.5}}
\put(2,0){\line(0,1){7}}
\put(0,-0.1){\makebox(0,0)[t]{$\ss 1$}}
\put(2,-0.1){\makebox(0,0)[t]{$\ss k$}}
\put(7.5,-0.1){\makebox(0,0)[t]{$\rightarrow i$}}
\put(-0.1,7.5){\makebox(0,0)[r]{$j\uparrow$}}
\put(5.6,5.6){\makebox(0,0)[lb]{$\ss i=j$}}
\put(2.5,0){\line(-1,1){0.5}}
\put(3,0){\line(-1,1){1}}
\put(3.5,0){\line(-1,1){1.5}}
\put(4,0){\line(-1,1){2}}
\put(4.5,0){\line(-1,1){2.25}}
\put(5,0){\line(-1,1){2.5}}
\put(5.5,0){\line(-1,1){2.75}}
\put(6,0){\line(-1,1){3}}
\put(6,0.5){\line(-1,1){2.75}}
\put(6,1){\line(-1,1){2.5}}
\put(6,1.5){\line(-1,1){2.25}}
\put(6,2){\line(-1,1){2}}
\end{picture}
\end{center}
$$\sum_{i=k}^{\infty}\;\sum_{j=1}^i=
\sum_{j=1}^{\infty}\;\sum_{i=\max(j,k)}^{\infty}$$
to find
$$\sum_{i=k}^{\infty}\frac{i!}{(i-k)!}a_i(x)=\sum_{j=1}^{\infty}(-1)^{j+1}
{\a+1 \ch j-1}x^j\sum_{i=\max(j,k)}^{\infty}\frac{(-1)^i}{(i-k)!}
{\a+2 \ch i-j}(\a+3)_{i-j}.$$
For $j=1,2,\ldots,k$ we obtain
\bea\sum_{i=k}^{\infty}\frac{(-1)^i}{(i-k)!}{\a+2 \ch i-j}(\a+3)_{i-j}
&=&\sum_{i=0}^{\infty}\frac{(-1)^{i+k}}{i!}{\a+2 \ch i+k-j}
(\a+3)_{i+k-j}\nn
&=&(-1)^j\sum_{i=0}^{\infty}\frac{(-\a-2)_{i+k-j}(\a+3)_{i+k-j}}
{i!(i+k-j)!}\n\eea
and for $j=k+1,k+2,\ldots$ we find by using the summation formula (\ref{sum})
\bea\sum_{i=j}^{\infty}\frac{(-1)^i}{(i-k)!}{\a+2 \ch i-j}(\a+3)_{i-j}
&=&\sum_{i=0}^{\infty}\frac{(-1)^{i+j}}{(i+j-k)!}{\a+2 \ch i}
(\a+3)_i\nn
&=&(-1)^j\sum_{i=0}^{\infty}\frac{(-\a-2)_i(\a+3)_i}{i!(i+j-k)!}\nn
&=&\frac{(-1)^j}{(j-k)!}\hyp{2}{1}{-\a-2,\a+3}{j-k+1}{1}\nn
&=&(-1)^j\frac{\G(j-k)}{\G(j-k+\a+3)\G(j-k-\a-2)}.\n\eea
Hence
\bea\sum_{i=k}^{\infty}\frac{i!}{(i-k)!}a_i(x)
&=&-\sum_{j=1}^k{\a+1 \ch j-1}x^j\sum_{i=0}^{\infty}
\frac{(-\a-2)_{i+k-j}(\a+3)_{i+k-j}}{i!(i+k-j)!}+{}\nn
& &{}-\sum_{j=k+1}^{\infty}{\a+1 \ch j-1}x^j\frac{\G(j-k)}
{\G(j-k+\a+3)\G(j-k-\a-2)}.\n\eea
For the last sum we find
\bea & &\sum_{j=k+1}^{\infty}{\a+1 \ch j-1}x^j\frac{\G(j-k)}
{\G(j-k+\a+3)\G(j-k-\a-2)}\nn
&=&\sum_{j=0}^{\infty}{\a+1 \ch j+k}x^{j+k+1}\frac{\G(j+1)}
{\G(j+\a+4)\G(j-\a-1)}\nn
&=&(-1)^kx^{k+1}\sum_{j=0}^{\infty}\frac{(-\a-1)_{j+k}}{(j+k)!}\frac
{\G(j+1)}{\G(j+\a+4)\G(j-\a-1)}(-x)^j\n\eea
which equals zero since
$$\frac{1}{\G(j-\a-1)}=0\;\mbox{ for }\;j=0,1,2,\ldots,\a+1$$
and
$$(-\a-1)_{j+k}=0\;\mbox{ for }\;j=\a+2,\a+3,\ldots\;
\mbox{ and }\;k=1,2,3,\ldots.$$
So we have
$$\sum_{i=k}^{\infty}\frac{i!}{(i-k)!}a_i(x)
=-\sum_{j=1}^k{\a+1 \ch j-1}x^j\sum_{i=0}^{\infty}
\frac{(-\a-2)_{i+k-j}(\a+3)_{i+k-j}}{i!(i+k-j)!}.$$
The inner sum equals zero if $k-j>\a+2$. Hence
\bea\sum_{i=k}^{\infty}\frac{i!}{(i-k)!}a_i(x)
&=&-\sum_{j=\max(1,k-\a-2)}^k{\a+1 \ch j-1}x^j
\frac{(-\a-2)_{k-j}(\a+3)_{k-j}}{(k-j)!}\times{}\nn
& &\hspace{4cm}\times\hyp{2}{1}{-\a-2+k-j,\a+3+k-j}{k-j+1}{1}\nn
&=&(-1)^{k+1}\sum_{j=\max(1,k-\a-2)}^k{\a+1 \ch j-1}
{\a+2 \ch k-j}(\a+3)_{k-j}(-x)^j\times{}\nn
& &\hspace{4cm}\times\hyp{2}{1}{-\a-2+k-j,\a+3+k-j}{k-j+1}{1}.\n\eea
Now we use the Vandermonde summation formula
$$\hyp{2}{1}{-n,b}{c}{1}=\frac{(c-b)_n}{(c)_n},\;n=0,1,2,\ldots$$
to find
$$\hyp{2}{1}{-\a-2+k-j,\a+3+k-j}{k-j+1}{1}
=\frac{(-\a-2)_{\a+2-k+j}}{(k-j+1)_{\a+2-k+j}}
=(-1)^{\a+2-k+j}.$$
Hence
\bea\sum_{i=k}^{\infty}\frac{i!}{(i-k)!}a_i(x)
&=&(-1)^{\a+1}\sum_{j=\max(1,k-\a-2)}^k{\a+1 \ch j-1}
{\a+2 \ch k-j}(\a+3)_{k-j}x^j\nn
&=&(-1)^{\a+1}\sum_{j=1}^k{\a+1 \ch j-1}{\a+2 \ch k-j}
(\a+3)_{k-j}x^j\nn
&=&(-1)^{\a+k}k!a_k(-x).\n\eea
This proves (\ref{suma2}).

\section{Some preliminary results for the Sobolev Laguerre polynomials}

In this section we look for a differential equation of the form
\bea\la{Sobdv} & &M\sum_{i=0}^{\infty}a_i(x)y^{(i)}(x)+
N\sum_{i=0}^{\infty}b_i(x)y^{(i)}(x)+{}\nn
& &\hspace{1cm}{}+MN\sum_{i=0}^{\infty}c_i(x)y^{(i)}(x)+
xy''(x)+(\a+1-x)y'(x)+ny(x)=0\eea
for the polynomials
$$y(x)=\gL=A_0\L+A_1\frac{d}{dx}\L+A_2\frac{d^2}{dx^2}\L,$$
where the coefficients $A_0$, $A_1$ and $A_2$ are defined by
\be\la{AAA}\l\{\ba{l}\ds A_0=1+M{n+\a \ch n-1}+\frac{n(\a+2)-(\a+1)}{(\a+1)(\a+3)}
N{n+\a \ch n-2}+{}\\
\ds\hspace{5cm}{}+\frac{MN}{(\a+1)(\a+2)}{n+\a \ch n-1}{n+\a+1 \ch n-2}\\  \\
\ds A_1=M{n+\a \ch n}+\frac{(n-1)}{(\a+1)}N{n+\a \ch n-1}+
\frac{2MN}{(\a+1)^2}{n+\a \ch n}{n+\a+1 \ch n-2}\\  \\
\ds A_2=\frac{N}{(\a+1)}{n+\a \ch n-1}+\frac{MN}{(\a+1)^2}{n+\a \ch n}
{n+\a+1 \ch n-1}.\ea\r.\ee
For details concerning these generalized Laguerre polynomials and their
definition the reader is referred to \cite{SIAM} and \cite{Thesis}.

Of course, since $L_n^{\a,M,0}(x)=\KL$, the coefficients $\l\{a_i(x)\r\}_{i=0}^{\infty}$
are given by (\ref{coeff}).

Although the general form is still an open problem so far, we know that the
differential equation given by (\ref{Sobdv}) is not unique as in the case of
the differential equation (\ref{dv}). We introduce the notation
$$\l\{\ba{l}b_0(0,\a,x)=0\\ \\
b_0(n,\a,x)=b_0(1,\a,x)+\b_0(n,\a,x),\;n=1,2,3,\ldots\\ \\
b_i(\a,x)=b_0(1,\a,x)b_i^*(\a,x)+\b_i(\a,x),\;i=1,2,3,\ldots\ea\r.$$
and
$$\l\{\ba{l}c_0(0,\a,x)=0\\ \\
c_0(n,\a,x)=b_0(1,\a,x)+\c_0(n,\a,x),\;n=1,2,3,\ldots\\ \\
c_i(\a,x)=b_0(1,\a,x)c_i^*(\a,x)+\c_i(\a,x),\;i=1,2,3,\ldots.\ea\r.$$
Now we can prove the following theorem~:

\st{The polynomials $\l\{\gL\r\}_{n=1}^{\infty}$ satisfy the following
infinite order differential equation~:
\be\la{dvstar}\sum_{i=0}^{\infty}b_i^*(\a,x)y^{(i)}(x)+
M\sum_{i=0}^{\infty}c_i^*(\a,x)y^{(i)}(x)=0,\ee
where
$$b_i^*(\a,x)=\frac{1}{i!}\sum_{j=0}^i(-1)^j{i \ch j}
(\a+1)_{i-j}x^j,\;i=0,1,2,\ldots$$
and
$$c_i^*(\a,x)=\frac{(-1)^i}{i!}x^i,\;i=0,1,2,\ldots.$$}

{\bf Proof.} The proof is very easy and is based on the observation that
\be\la{sumbstar}\sum_{i=0}^{\infty}b_i^*(\a,x)D^{i+k}\L
=\frac{(-n)_k}{n\G(k)},\;n\ge 1,\;k=0,1,2\ee
and
\be\la{sumcstar}\sum_{i=0}^{\infty}c_i^*(\a,x)D^{i+k}\L
={n+\a \ch n}\frac{(-n)_k}{(\a+1)_k},\;k=0,1,2.\ee

To prove (\ref{sumbstar}) we change the order of summation to obtain
\bea\sum_{i=0}^{\infty}b_i^*(\a,x)D^{i+k}\L
&=&\sum_{i=0}^{\infty}\sum_{j=0}^i\frac{(-1)^j}{i!}{i \ch j}(\a+1)_{i-j}x^jD^{i+k}\L\nn
&=&\sum_{j=0}^{\infty}\sum_{i=j}^{\infty}\frac{(-1)^j}{i!}{i \ch j}
(\a+1)_{i-j}x^jD^{i+k}\L\nn
&=&\sum_{j=0}^{\infty}\sum_{i=0}^{\infty}\frac{(-1)^j}{(i+j)!}{i+j \ch j}
(\a+1)_ix^jD^{i+j+k}\L\nn
&=&\sum_{i=0}^{\infty}\frac{(\a+1)_i}{i!}\sum_{j=0}^{\infty}\frac{(-1)^j}{j!}
x^jD^{i+j+k}\L.\n\eea
Now we use the definition of the classical Laguerre polynomial
$$\L={n+\a \ch n}\hyp{1}{1}{-n}{\a+1}{x}$$
to obtain
\bea\sum_{j=0}^{\infty}\frac{(-1)^j}{j!}x^jD^{i+j+k}\L
&=&{n+\a \ch n}\sum_{j=0}^{\infty}\frac{(-1)^j}{j!}x^j\sum_{m=i+j+k}^{\infty}
\frac{(-n)_m}{(\a+1)_m}\frac{x^{m-i-j-k}}{(m-i-j-k)!}\nn
&=&{n+\a \ch n}\sum_{j=0}^{\infty}\sum_{m=j}^{\infty}\frac{(-1)^j}{j!}
\frac{(-n)_{m+i+k}}{(\a+1)_{m+i+k}}\frac{x^m}{(m-j)!}\nn
&=&{n+\a \ch n}\sum_{m=0}^{\infty}\frac{(-n)_{m+i+k}}{(\a+1)_{m+i+k}}
\frac{x^m}{m!}\sum_{j=0}^m(-1)^j{m \ch j}\nn
&=&{n+\a \ch n}\frac{(-n)_{i+k}}{(\a+1)_{i+k}}.\n\eea
Hence, by using the summation formula (\ref{sum}) we find
\bea\sum_{i=0}^{\infty}b_i^*(\a,x)D^{i+k}\L
&=&{n+\a \ch n}\sum_{i=0}^{\infty}\frac{(\a+1)_i}{i!}\frac{(-n)_{i+k}}{(\a+1)_{i+k}}\nn
&=&{n+\a \ch n}\frac{(-n)_k}{(\a+1)_k}\hyp{2}{1}{-n+k,\a+1}{\a+k+1}{1}\nn
&=&{n+\a \ch n}\frac{(-n)_k}{(\a+1)_k}\frac{\G(n)\G(\a+k+1)}{\G(n+\a+1)\G(k)}
=\frac{(-n)_k}{n\G(k)},\n\eea
which proves (\ref{sumbstar}). The proof of (\ref{sumcstar}) is much shorter~:
\bea\sum_{i=0}^{\infty}c_i^*(\a,x)D^{i+k}\L
&=&{n+\a \ch n}\sum_{i=0}^{\infty}\frac{(-1)^i}{i!}x^i\sum_{m=i+k}^{\infty}
\frac{(-n)_m}{(\a+1)_m}\frac{x^{m-i-k}}{(m-i-k)!}\nn
&=&{n+\a \ch n}\sum_{i=0}^{\infty}\frac{(-1)^i}{i!}x^i\sum_{m=i}^{\infty}
\frac{(-n)_{m+k}}{(\a+1)_{m+k}}\frac{x^{m-i}}{(m-i)!}\nn
&=&{n+\a \ch n}\sum_{m=0}^{\infty}\sum_{i=0}^m\frac{(-1)^i}{i!}
\frac{(-n)_{m+k}}{(\a+1)_{m+k}}\frac{x^m}{(m-i)!}\nn
&=&{n+\a \ch n}\sum_{m=0}^{\infty}\frac{(-n)_{m+k}}{(\a+1)_{m+k}}
\frac{x^m}{m!}\sum_{i=0}^m(-1)^i{m \ch i}\nn
&=&{n+\a \ch n}\frac{(-n)_k}{(\a+1)_k}.\n\eea

By using these formulas (\ref{sumbstar}) and (\ref{sumcstar}) and by using
the definition (\ref{AAA}) of the coefficients $A_0$, $A_1$ and $A_2$ we find
\bea & &\!\!\!\!\!A_0\sum_{i=0}^{\infty}b_i^*(\a,x)D^i\L+
A_1\sum_{i=0}^{\infty}b_i^*(\a,x)D^{i+1}\L
+A_2\sum_{i=0}^{\infty}b_i^*(\a,x)D^{i+2}\L+{}\nn
& &\!\!\!\!\!{}+MA_0\sum_{i=0}^{\infty}c_i^*(\a,x)D^i\L+
MA_1\sum_{i=0}^{\infty}c_i^*(\a,x)D^{i+1}\L+{}\nn
& &\hspace{7cm}{}+MA_2\sum_{i=0}^{\infty}c_i^*(\a,x)D^{i+2}\L\nn
&=&-A_1+(n-1)A_2+M\l[{n+\a \ch n}A_0-{n+\a \ch n-1}A_1+{n+\a \ch n-2}A_2\r]=0.\n\eea
This proves (\ref{dvstar}).

\vspace{5mm}

In the special cases $\a=0$, $\a=1$ and $\a=2$ differential equations of the
form (\ref{Sobdv}) are found too. In these
three special cases of integer values of the parameter $\a$ we find a linear
differential equation of formal order $4\a+10$. By formal order we mean that
for special cases ($M=0$ or $N=0$) the true order might be lower.
We give the results, but we will not give any proofs here.

\subsection{The special case $\a=0$}

If we take $\a=0$ in (\ref{coeff}) we find
$$a_0(n,0,x)=\frac{1}{2}n(n+1),\;\ndots$$
and
\bea a_1(0,x)&=&-x\nn
a_2(0,x)&=&3x-\frac{1}{2}x^2\nn
a_3(0,x)&=&-2x+x^2\nn
a_4(0,x)&=&-\frac{1}{2}x^2\nn
a_i(0,x)&=&0,\;i=5,6,7,\ldots.\n\eea
For the coefficients $\l\{\b_i(x)\r\}_{i=0}^{\infty}$ we find in this case
$$\b_0(n,0,x)=\frac{1}{12}n^2(n^2-1),\;n=1,2,3,\ldots$$
and
\bea \b_1(0,x)&=&0\nn
\b_2(0,x)&=&1-\frac{1}{2}x^2\nn
\b_3(0,x)&=&-3-x+\frac{9}{2}x^2-\frac{1}{2}x^3\nn
\b_4(0,x)&=&2+3x-\frac{25}{2}x^2+\frac{17}{6}x^3-\frac{1}{12}x^4\nn
\b_5(0,x)&=&-2x+\frac{27}{2}x^2-\frac{11}{2}x^3+\frac{1}{3}x^4\nn
\b_6(0,x)&=&-5x^2+\frac{9}{2}x^3-\frac{1}{2}x^4\nn
\b_7(0,x)&=&-\frac{4}{3}x^3+\frac{1}{3}x^4\nn
\b_8(0,x)&=&-\frac{1}{12}x^4\nn
\b_i(0,x)&=&0,\;i=9,10,11,\ldots.\n\eea
Finally, the coefficients $\l\{\c_i(x)\r\}_{i=0}^{\infty}$ turn out to be
$$\c_0(n,0,x)=\frac{1}{120}n(n^2-1)(n+2)(2n+1),\;n=1,2,3,\ldots$$
and
\bea \c_1(0,x)&=&0\nn
\c_2(0,x)&=&-\frac{1}{2}x^2\nn
\c_3(0,x)&=&5x^2-\frac{2}{3}x^3\nn
\c_4(0,x)&=&-\frac{35}{2}x^2+5x^3-\frac{5}{24}x^4\nn
\c_5(0,x)&=&28x^2-14x^3+\frac{5}{4}x^4-\frac{1}{60}x^5\nn
\c_6(0,x)&=&-21x^2+\frac{56}{3}x^3-\frac{35}{12}x^4+\frac{1}{12}x^5\nn
\c_7(0,x)&=&6x^2-12x^3+\frac{10}{3}x^4-\frac{1}{6}x^5\nn
\c_8(0,x)&=&3x^3-\frac{15}{8}x^4+\frac{1}{6}x^5\nn
\c_9(0,x)&=&\frac{5}{12}x^4-\frac{1}{12}x^5\nn
\c_{10}(0,x)&=&\frac{1}{60}x^5\nn
\c_i(0,x)&=&0,\;i=11,12,13,\ldots.\n\eea

Note that we have in this case~:
$$\sum_{i=1}^4a_i(0,x)=\sum_{i=1}^8\b_i(0,x)=\sum_{i=1}^{10}\c_i(0,x)=0.$$

Hence, we have found the following linear differential equation of formal
order $10$~:
\bea & &\frac{1}{60}MNx^5y^{(10)}(x)+\frac{1}{12}MN(5x^4-x^5)y^{(9)}(x)+{}\nn
& &{}+\l[\frac{1}{24}MN(72x^3-45x^4+4x^5)-\frac{1}{12}Nx^4\r]y^{(8)}(x)+{}\nn
& &{}+\l[\frac{1}{6}MN(36x^2-72x^3+20x^4-x^5)-\frac{1}{3}N(4x^3-x^4)\r]
y^{(7)}(x)+{}\nn
& &{}+\l[\frac{1}{12}MN(-252x^2+224x^3-35x^4+x^5)
+\frac{1}{2}N(-10x^2+9x^3-x^4)\r]y^{(6)}(x)+{}\nn
& &{}+\l[\frac{1}{60}MN(1680x^2-840x^3+75x^4-x^5)
+\frac{1}{6}N(-12x+81x^2-33x^3+2x^4)\r]y^{(5)}(x)+{}\nn
& &{}+\l[\frac{1}{24}MN(-420x^2+120x^3-5x^4)+{}\r.\nn
& &\hspace{3cm}\l.{}+\frac{1}{12}N(24+36x-150x^2+34x^3-x^4)
-\frac{1}{2}Mx^2\r]y^{(4)}(x)+{}\nn
& &{}+\l[\frac{1}{3}MN(15x^2-2x^3)
+\frac{1}{2}N(-6-2x+9x^2-x^3)+M(-2x+x^2)\r]y^{(3)}(x)+{}\nn
& &{}+\l[-\frac{1}{2}MNx^2+\frac{1}{2}N(2-x^2)+\frac{1}{2}M(6x-x^2)+x\r]y''(x)
+\l[1-(M+1)x\r]y'(x)+{}\nn
& &+\frac{1}{120}n\l[MN(n^2-1)(n+2)(2n+1)+10Nn(n^2-1)+60M(n+1)+120\r]y(x)=0,\n\eea
for the polynomials
$$y(x):=L_n^{0,M,N}(x)=A_0L_n(x)+A_1L_n'(x)+A_2L_n''(x)$$
with
$$L_n(x):=L_n^{(0)}(x)=\sum_{k=0}^n\frac{(-1)^k}{k!}{n \ch k}x^k,\;\ndots$$
and
$$\l\{\ba{l}\ds A_0=1+Mn+\frac{1}{6}Nn(n-1)(2n-1)+\frac{1}{12}MNn^2(n^2-1)\\ \\
\ds A_1=M+Nn(n-1)+\frac{1}{3}MNn(n^2-1)\\ \\
\ds A_2=Nn+\frac{1}{2}MNn(n+1).\ea\r.$$

\subsection{The special case $\a=1$}

In the special case that $\a=1$ we find successively
$$a_0(n,1,x)=\frac{1}{6}n(n+1)(n+2),\;\ndots,$$
\bea a_1(1,x)&=&-x\nn
a_2(1,x)&=&6x-x^2\nn
a_3(1,x)&=&-10x+4x^2-\frac{1}{6}x^3\nn
a_4(1,x)&=&5x-5x^2+\frac{1}{2}x^3\nn
a_5(1,x)&=&2x^2-\frac{1}{2}x^3\nn
a_6(1,x)&=&\frac{1}{6}x^3\nn
a_i(1,x)&=&0,\;i=7,8,9,\ldots,\n\eea
$$\b_0(n,1,x)=\frac{1}{240}n(n^2-1)(n+2)(3n-1),\;n=1,2,3,\ldots,$$
\bea \b_1(1,x)&=&0\nn
\b_2(1,x)&=&\frac{3}{2}-\frac{1}{4}x^2\nn
\b_3(1,x)&=&-8-\frac{3}{2}x+4x^2-\frac{5}{12}x^3\nn
\b_4(1,x)&=&\frac{25}{2}+8x-\frac{77}{4}x^2+\frac{47}{12}x^3-\frac{7}{48}x^4\nn
\b_5(1,x)&=&-6-\frac{25}{2}x+\frac{77}{2}x^2-\frac{51}{4}x^3
+\frac{23}{24}x^4-\frac{1}{80}x^5\nn
\b_6(1,x)&=&6x-34x^2+\frac{227}{12}x^3-\frac{19}{8}x^4+\frac{1}{16}x^5\nn
\b_7(1,x)&=&11x^2-\frac{79}{6}x^3+\frac{17}{6}x^4-\frac{1}{8}x^5\nn
\b_8(1,x)&=&\frac{7}{2}x^3-\frac{79}{48}x^4+\frac{1}{8}x^5\nn
\b_9(1,x)&=&\frac{3}{8}x^4-\frac{1}{16}x^5\nn
\b_{10}(1,x)&=&\frac{1}{80}x^5\nn
\b_i(1,x)&=&0,\;i=11,12,13,\ldots,\n\eea
$$\c_0(n,1,x)=\frac{n(n^2-1)(n+2)(n+3)(5n^2+10n+2)}{10080},\;n=1,2,3,\ldots$$
and
\bea \c_1(1,x)&=&0\nn
\c_2(1,x)&=&-\frac{1}{4}x^2\nn
\c_3(1,x)&=&5x^2-\frac{2}{3}x^3\nn
\c_4(1,x)&=&-35x^2+\frac{115}{12}x^3-\frac{23}{48}x^4\nn
\c_5(1,x)&=&119x^2-\frac{105}{2}x^3+\frac{65}{12}x^4-\frac{31}{240}x^5\nn
\c_6(1,x)&=&-\frac{441}{2}x^2+147x^3-\frac{49}{2}x^4+\frac{29}{24}x^5
-\frac{1}{72}x^6\nn
\c_7(1,x)&=&228x^2-232x^3+\frac{117}{2}x^4-\frac{14}{3}x^5
+\frac{1}{9}x^6-\frac{1}{2016}x^7\nn
\c_8(1,x)&=&-\frac{495}{4}x^2+\frac{837}{4}x^3-\frac{1289}{16}x^4
+\frac{467}{48}x^5-\frac{3}{8}x^6+\frac{1}{288}x^7\nn
\c_9(1,x)&=&\frac{55}{2}x^2-\frac{605}{6}x^3+\frac{129}{2}x^4
-\frac{571}{48}x^5+\frac{25}{36}x^6-\frac{1}{96}x^7\nn
\c_{10}(1,x)&=&\frac{121}{6}x^3-\frac{671}{24}x^4+\frac{257}{30}x^5
-\frac{55}{72}x^6+\frac{5}{288}x^7\nn
\c_{11}(1,x)&=&\frac{61}{12}x^4-\frac{27}{8}x^5+\frac{1}{2}x^6-\frac{5}{288}x^7\nn
\c_{12}(1,x)&=&\frac{9}{16}x^5-\frac{13}{72}x^6+\frac{1}{96}x^7\nn
\c_{13}(1,x)&=&\frac{1}{36}x^6-\frac{1}{288}x^7\nn
\c_{14}(1,x)&=&\frac{1}{2016}x^7\nn
\c_i(1,x)&=&0,\;i=15,16,17,\ldots.\n\eea

Hence, in this case we have
$$\sum_{i=1}^6a_i(1,x)=\sum_{i=1}^{10}\b_i(1,x)=\sum_{i=1}^{14}\c_i(1,x)=0.$$

This implies that we have found a linear differential equation of formal order
$14$ for the polynomials $\set{L_n^{1,M,N}(x)}$.

\subsection{The special case $\a=2$}

If we take $\a=2$ we find successively
$$a_0(n,2,x)=\frac{1}{24}n(n+1)(n+2)(n+3),\;\ndots,$$
\bea a_1(2,x)&=&-x\nn
a_2(2,x)&=&10x-\frac{3}{2}x^2\nn
a_3(2,x)&=&-30x+10x^2-\frac{1}{2}x^3\nn
a_4(2,x)&=&35x-\frac{45}{2}x^2+\frac{5}{2}x^3-\frac{1}{24}x^4\nn
a_5(2,x)&=&-14x+21x^2-\frac{9}{2}x^3+\frac{1}{6}x^4\nn
a_6(2,x)&=&-7x^2+\frac{7}{2}x^3-\frac{1}{4}x^4\nn
a_7(2,x)&=&-x^3+\frac{1}{6}x^4\nn
a_8(2,x)&=&-\frac{1}{24}x^4\nn
a_i(2,x)&=&0,\;i=9,10,11,\ldots,\n\eea
$$\b_0(n,2,x)=\frac{1}{1080}n(n^2-1)(n+2)(n+3)(2n-1),\;n=1,2,3,\ldots,$$
\bea \b_1(2,x)&=&0\nn
\b_2(2,x)&=&2-\frac{1}{6}x^2\nn
\b_3(2,x)&=&-\frac{50}{3}-2x+\frac{25}{6}x^2-\frac{7}{18}x^3\nn
\b_4(2,x)&=&45+\frac{50}{3}x-\frac{92}{3}x^2+\frac{11}{2}x^3-\frac{5}{24}x^4\nn
\b_5(2,x)&=&-49-45x+\frac{290}{3}x^2-27x^3+\frac{49}{24}x^4-\frac{13}{360}x^5\nn
\b_6(2,x)&=&\frac{56}{3}+49x-\frac{887}{6}x^2+\frac{188}{3}x^3-\frac{91}{12}x^4
+\frac{97}{360}x^5-\frac{1}{540}x^6\nn
\b_7(2,x)&=&-\frac{56}{3}x+\frac{217}{2}x^2-\frac{451}{6}x^3+\frac{169}{12}x^4
-\frac{29}{36}x^5+\frac{1}{90}x^6\nn
\b_8(2,x)&=&-\frac{92}{3}x^2+\frac{271}{6}x^3-\frac{337}{24}x^4
+\frac{5}{4}x^5-\frac{1}{36}x^6\nn
\b_9(2,x)&=&-\frac{97}{9}x^3+\frac{173}{24}x^4-\frac{77}{72}x^5+\frac{1}{27}x^6\nn
\b_{10}(2,x)&=&-\frac{3}{2}x^4+\frac{173}{360}x^5-\frac{1}{36}x^6\nn
\b_{11}(2,x)&=&-\frac{4}{45}x^5+\frac{1}{90}x^6\nn
\b_{12}(2,x)&=&-\frac{1}{540}x^6\nn
\b_i(2,x)&=&0,\;i=13,14,15,\ldots,\n\eea
$$\c_0(n,2,x)=\frac{n(n^2-1)(n+2)(n+3)(n+4)(2n+3)(7n^2+21n+2)}{1088640},
\;n=1,2,3,\ldots$$
and
\bea \c_1(2,x)&=&0\nn
\c_2(2,x)&=&-\frac{1}{6}x^2\nn
\c_3(2,x)&=&\frac{35}{6}x^2-\frac{13}{18}x^3\nn
\c_4(2,x)&=&-70x^2+\frac{35}{2}x^3-\frac{7}{8}x^4\nn
\c_5(2,x)&=&413x^2-161x^3+\frac{49}{3}x^4-\frac{13}{30}x^5\nn
\c_6(2,x)&=&-1386x^2+\frac{2317}{3}x^3-\frac{245}{2}x^4+\frac{119}{18}x^5
-\frac{73}{720}x^6\nn
\c_7(2,x)&=&2827x^2-2189x^3+497x^4-42x^5+\frac{21}{16}x^6-\frac{59}{5040}x^7\nn
\c_8(2,x)&=&-3575x^2+3872x^3-\frac{9779}{8}x^4+148x^5-\frac{29}{4}x^6
+\frac{19}{144}x^7-\frac{11}{17280}x^8\nn
\c_9(2,x)&=&\frac{16445}{6}x^2-\frac{77935}{18}x^3+\frac{11473}{6}x^4
-\frac{5797}{18}x^5+\frac{3257}{144}x^6-\frac{31}{48}x^7+{}\nn
& &\hspace{8cm}{}+\frac{11}{1728}x^8
-\frac{1}{77760}x^9\nn
\c_{10}(2,x)&=&-\frac{7007}{6}x^2+\frac{17875}{6}x^3-\frac{23023}{12}x^4
+\frac{4521}{10}x^5-\frac{2651}{60}x^6+\frac{1303}{720}x^7+{}\nn
& &\hspace{8cm}{}-\frac{121}{4320}x^8
+\frac{1}{8640}x^9\nn
\c_{11}(2,x)&=&\frac{637}{3}x^2-\frac{3458}{3}x^3+\frac{7189}{6}x^4
-\frac{1235}{3}x^5+\frac{13511}{240}x^6-\frac{2311}{720}x^7+{}\nn
& &\hspace{8cm}{}+\frac{77}{1080}x^8
-\frac{1}{2160}x^9\nn
\c_{12}(2,x)&=&\frac{1729}{9}x^3-\frac{1274}{3}x^4+\frac{2119}{9}x^5
-\frac{845}{18}x^6+\frac{899}{240}x^7-\frac{1001}{8640}x^8+\frac{7}{6480}x^9\nn
\c_{13}(2,x)&=&\frac{196}{3}x^4-77x^5+\frac{1189}{48}x^6-\frac{415}{144}x^7
+\frac{539}{4320}x^8-\frac{7}{4320}x^9\nn
\c_{14}(2,x)&=&11x^5-\frac{361}{48}x^6+\frac{1423}{1008}x^7
-\frac{77}{864}x^8+\frac{7}{4320}x^9\nn
\c_{15}(2,x)&=&\frac{361}{360}x^6-\frac{2}{5}x^7+\frac{11}{270}x^8
-\frac{7}{6480}x^9\nn
\c_{16}(2,x)&=&\frac{1}{20}x^7-\frac{187}{17280}x^8+\frac{1}{2160}x^9\nn
\c_{17}(2,x)&=&\frac{11}{8640}x^8-\frac{1}{8640}x^9\nn
\c_{18}(2,x)&=&\frac{1}{77760}x^9\nn
\c_i(2,x)&=&0,\;i=19,20,21,\ldots.\n\eea

Hence, we have
$$\sum_{i=1}^8a_i(2,x)=\sum_{i=1}^{12}\b_i(2,x)=\sum_{i=1}^{18}\c_i(2,x)=0.$$

This implies that the polynomials $\set{L_n^{2,M,N}(x)}$ satisfy a linear
differential equation of formal order $18$.

\section{The generalized Jacobi polynomials}

In this section we will deal with the problem of finding a differential
equation for the generalized Jacobi polynomials $\set{\GP}$.

Since the well-known second order differential equation for the classical
Jacobi polynomials $\set{\P}$ is given by
$$(1-x^2)y''(x)+\l[\b-\a-(\a+\b+2)x\r]y'(x)+n(n+\a+\b+1)y(x)=0,$$
it is clear that we look for a differential equation of the form
\bea & &M\sum_{i=0}^{\infty}a_i(x)y^{(i)}(x)+N\sum_{i=0}^{\infty}b_i(x)y^{(i)}(x)
+MN\sum_{i=0}^{\infty}c_i(x)y^{(i)}(x)+{}\nn
& &\hspace{1cm}{}+(1-x^2)y''(x)
+\l[\b-\a-(\a+\b+2)x\r]y'(x)+n(n+\a+\b+1)y(x)=0,\n\eea
where
$$y(x)=P_n^{\a,\b,M,N}(x)=A_0P_n^{(\a,\b)}(x)+\l[A_1(1-x)-A_2(1+x)\r]\frac{d}{dx}P_n^{(\a,\b)}(x)$$
and
$$\l\{\ba{l}\ds A_0=1+M\frac{\ds {n+\b \ch n-1}{n+\a+\b+1 \ch n}}{\ds {n+\a \ch n}}
+N\frac{\ds {n+\a \ch n-1}{n+\a+\b+1 \ch n}}{\ds {n+\b \ch n}}+{}\\  \\
\ds\hspace{7cm}{}+MN\frac{(\a+\b+2)^2}{(\a+1)(\b+1)}{n+\a+\b+1 \ch n-1}^2\\ \\
\ds A_1=\frac{M}{(\a+\b+1)}\frac{\ds {n+\b \ch n}{n+\a+\b \ch n}}{\ds {n+\a \ch n}}
+\frac{MN}{(\a+1)}{n+\a+\b \ch n-1}{n+\a+\b+1 \ch n}\\ \\
\ds A_2=\frac{N}{(\a+\b+1)}\frac{\ds {n+\a \ch n}{n+\a+\b \ch n}}{\ds {n+\b \ch n}}
+\frac{MN}{(\b+1)}{n+\a+\b \ch n-1}{n+\a+\b+1 \ch n}.\ea\r.$$
Here we used the same definition as in \cite{Koorn}, but in a slightly different
notation.

In this case the differential equation will be unique in its general form if
it exists. We introduce the notation
$$\l\{\ba{l}a_0(x)=a_0(n,\a,\b,x)\\ \\
a_i(x)=a_i(\a,\b,x),\;i=1,2,3,\ldots,\ea\r.$$
$$\l\{\ba{l}b_0(x)=b_0(n,\a,\b,x)\\ \\
b_i(x)=b_i(\a,\b,x),\;i=1,2,3,\ldots\ea\r.$$
and
$$\l\{\ba{l}c_0(x)=c_0(n,\a,\b,x)\\ \\
c_i(x)=c_i(\a,\b,x),\;i=1,2,3,\ldots.\ea\r.$$
Since the polynomials $\set{\GP}$ satisfy the symmetry relation
$$P_n^{\a,\b,M,N}(-x)=(-1)^nP_n^{\b,\a,N,M}(x)$$
we have
$$\l\{\ba{l}a_0(n,\a,\b,-x)=b_0(n,\b,\a,x)\\  \\
a_i(\a,\b,-x)=(-1)^ib_i(\b,\a,x),\;i=1,2,3,\ldots\ea\r.$$
and
$$\l\{\ba{l}c_0(n,\a,\b,-x)=c_0(n,\b,\a,x)\\  \\
c_i(\a,\b,-x)=(-1)^ic_i(\b,\a,x),\;i=1,2,3,\ldots.\ea\r.$$

The general form is still an open problem, but the symmetric case $\a=\b$
and $M=N$ seems to be much less difficult. All special cases known so far
(see for instance \cite{Conj}) point out that if $\a=\b$ and $M=N$ we can
choose $c_i(x)=0$ for all $i=0,1,2,\ldots$.
Therefore, we will consider this special case first.

\subsection{The symmetric generalized Jacobi polynomials}

All known examples seem to point out that there might be a differential
equation of the form
\be\la{dvSGP}M\sum_{i=0}^{\infty}a_i(x)y^{(i)}(x)+
(1-x^2)y''(x)-2(\a+1)xy'(x)+n(n+2\a+1)y(x)=0,\ee
where
\be\la{defSGP}y(x)=\SGP=C_0\SP-C_1x\frac{d}{dx}\SP\ee
and
\be\la{CC}\l\{\ba{l}\ds C_0=1+M\frac{2n}{(\a+1)}{n+2\a+1 \ch n}
+4M^2{n+2\a+1 \ch n-1}^2\\
\\
\ds C_1=\frac{2M}{(2\a+1)}{n+2\a \ch n}
+\frac{2M^2}{(\a+1)}{n+2\a \ch n-1}{n+2\a+1 \ch n}.\ea\r.\ee

This differential equation turns out not to be unique. We write
$$\l\{\ba{l}a_0(x):=a_0(n,\a,x),\;\ndots\\  \\
a_i(x):=a_i(\a,x),\;i=1,2,3,\ldots.\ea\r.$$

If we substitute (\ref{defSGP}) in the differential equation (\ref{dvSGP})
then we finally find three equations for the coefficients
$\l\{a_i(x)\r\}_{i=0}^{\infty}$ which are equivalent to the following two~:
$$\sum_{i=0}^{\infty}a_i(x)D^i\SP=\frac{4}{(2\a+1)}{n+2\a \ch n}
\frac{d^2}{dx^2}\SP$$
and
$$\sum_{i=0}^{\infty}ia_i(x)D^i\SP+x\sum_{i=0}^{\infty}a_i(x)D^{i+1}\SP
=4{n+2\a+1 \ch n-1}\frac{d^2}{dx^2}\SP.$$

We introduce the notation
$$\l\{\ba{l}a_0(n,\a,x)=a_0(1,\a,x)b_0(n,\a,x)+c_0(n,\a,x),\;\ndots\\ \\
a_i(\a,x)=a_0(1,\a,x)b_i(\a,x)+c_i(\a,x),\;i=1,2,3,\ldots.\ea\r.$$

Now we can prove the following theorem.

\st{The polynomials $\set{\SGP}$ satisfy a linear infinite order differential
equation of the form~:
$$\sum_{i=0}^{\infty}b_i(x)y^{(i)}(x)=0,$$
where
$$\l\{\ba{l}\ds b_0(x):=b_0(n,\a,x)=\frac{1}{2}\l[1-(-1)^n\r],\;\ndots\\
\\
\ds b_i(x):=b_i(\a,x)=\frac{2^{i-1}}{i!}(-x)^i,\;i=1,2,3,\ldots.\ea\r.$$}

{\bf Proof.} To prove this theorem we have to show that
$$\l\{\ba{l}\ds \sum_{i=0}^{\infty}b_i(x)D^i\SP=0\\ \\
\ds \sum_{i=0}^{\infty}ib_i(x)D^i\SP+x\sum_{i=0}^{\infty}b_i(x)D^{i+1}\SP=0.\ea\r.$$
To do this we note that the polynomials $\set{\SP}$ are defined by
$$\SP={n+\a \ch n}\hyp{2}{1}{-n,n+2\a+1}{\a+1}{\frac{1-x}{2}}.$$
Further we write
$$\sum_{i=0}^{\infty}b_i(x)D^i\SP=\frac{1}{2}\l[1-(-1)^n\r]\SP+\sum_{i=1}^{\infty}b_i(x)D^i\SP$$
and
$$\sum_{i=0}^{\infty}b_i(x)D^{i+1}\SP=\frac{1}{2}\l[1-(-1)^n\r]\frac{d}{dx}\SP
+\sum_{i=1}^{\infty}b_i(x)D^{i+1}\SP.$$
Now we easily obtain
\bea\sum_{i=1}^{\infty}b_i(x)D^i\SP
&=&\frac{1}{2}{n+\a \ch n}\sum_{i=1}^{\infty}\frac{x^i}{i!}\sum_{k=i}^{\infty}
\frac{(-n)_k(n+2\a+1)_k}{(\a+1)_k(k-i)!}\l(\frac{1-x}{2}\r)^{k-i}\nn
&=&\frac{1}{2}{n+\a \ch n}\sum_{k=1}^{\infty}\frac{(-n)_k(n+2\a+1)_k}{(\a+1)_kk!}
\sum_{i=1}^k{k \ch i}x^i\l(\frac{1-x}{2}\r)^{k-i}\nn
&=&\frac{1}{2}{n+\a \ch n}\sum_{k=1}^{\infty}\frac{(-n)_k(n+2\a+1)_k}{(\a+1)_kk!}
\l[\l(\frac{1+x}{2}\r)^k-\l(\frac{1-x}{2}\r)^k\r]\nn
&=&\frac{1}{2}\l[P_n^{(\a,\a)}(-x)-\SP\r]=-\frac{1}{2}\l[1-(-1)^n\r]\SP,\n\eea
\bea\sum_{i=0}^{\infty}ib_i(x)D^i\SP
&=&\frac{1}{2}{n+\a \ch n}\sum_{i=1}^{\infty}\frac{x^i}{(i-1)!}\sum_{k=i}^{\infty}
\frac{(-n)_k(n+2\a+1)_k}{(\a+1)_k(k-i)!}\l(\frac{1-x}{2}\r)^{k-i}\nn
&=&\frac{1}{2}{n+\a \ch n}\sum_{k=1}^{\infty}\frac{(-n)_k(n+2\a+1)_k}{(\a+1)_k(k-1)!}
\sum_{i=1}^k{k-1 \ch i-1}x^i\l(\frac{1-x}{2}\r)^{k-i}\nn
&=&\frac{1}{2}{n+\a \ch n}x\sum_{k=1}^{\infty}\frac{(-n)_k(n+2\a+1)_k}{(\a+1)_k(k-1)!}
\l(\frac{1+x}{2}\r)^{k-1}\nn
&=&(-1)^nx\frac{d}{dx}\SP\n\eea
and
\bea & &\sum_{i=1}^{\infty}b_i(x)D^{i+1}\SP\nn
&=&-\frac{1}{4}{n+\a \ch n}\sum_{i=1}^{\infty}\frac{x^i}{i!}\sum_{k=i}^{\infty}
\frac{(-n)_{k+1}(n+2\a+1)_{k+1}}{(\a+1)_{k+1}(k-i)!}\l(\frac{1-x}{2}\r)^{k-i}\nn
&=&-\frac{1}{4}{n+\a \ch n}\sum_{k=1}^{\infty}\frac{(-n)_{k+1}(n+2\a+1)_{k+1}}
{(\a+1)_{k+1}k!}\sum_{i=1}^k{k \ch i}x^i\l(\frac{1-x}{2}\r)^{k-i}\nn
&=&-\frac{1}{4}{n+\a \ch n}\sum_{k=1}^{\infty}\frac{(-n)_{k+1}(n+2\a+1)_{k+1}}
{(\a+1)_{k+1}k!}\l[\l(\frac{1+x}{2}\r)^k-\l(\frac{1-x}{2}\r)^k\r]\nn
&=&-\frac{1}{4}{n+\a \ch n}\sum_{k=0}^{\infty}\frac{(-n)_{k+1}(n+2\a+1)_{k+1}}
{(\a+1)_{k+1}k!}\l[\l(\frac{1+x}{2}\r)^k-\l(\frac{1-x}{2}\r)^k\r]\nn
&=&-\frac{1}{4}{n+\a \ch n}\sum_{k=1}^{\infty}\frac{(-n)_k(n+2\a+1)_k}
{(\a+1)_k(k-1)!}\l[\l(\frac{1+x}{2}\r)^{k-1}-\l(\frac{1-x}{2}\r)^{k-1}\r]\nn
&=&-\frac{1}{2}\l[1+(-1)^n\r]\frac{d}{dx}\SP.\n\eea
This proves the theorem.

For the coefficients $\l\{c_i(x)\r\}_{i=0}^{\infty}$ we find by using packages
like Maple or Mathematica~:
\bea c_0(0,\a,x)&=&0\nn
c_0(1,\a,x)&=&0\nn
c_0(2,\a,x)&=&4(2\a+3)\nn
c_0(3,\a,x)&=&4(2\a+3)(2\a+5)\nn
c_0(4,\a,x)&=&2(2\a+3)(2\a+5)(2\a+6)\nn
c_0(5,\a,x)&=&\frac{2}{3}(2\a+3)(2\a+5)(2\a+6)(2\a+7).\n\eea
Hence, we might guess that
\be\la{cnul}c_0(n,\a,x)=\frac{2n(n-1)}{(\a+2)}{n+2\a+2 \ch n}=4(2\a+3)
{n+2\a+2 \ch n-2},\;\ndots.\ee
If we write
\be\la{c}c_i(\a,x)=(2\a+3)(1-x^2)c_i^*(\a,x),\;i=1,2,3,\ldots\ee
we obtain
\bea c_1^*(\a,x)&=&0\nn
c_2^*(\a,x)&=&2\nn
c_3^*(\a,x)&=&\frac{4}{3}(\a+1)x\nn
c_4^*(\a,x)&=&\frac{1}{6}(\a+1)\l[(2\a+1)x^2-1\r]\nn
c_5^*(\a,x)&=&\frac{1}{45}\a(\a+1)x\l[(2\a+1)x^2-3\r]\nn
c_6^*(\a,x)&=&\frac{1}{1080}\a(\a+1)\l[(2\a-1)(2\a+1)x^4-6(2\a-1)x^2+3\r]\nn
c_7^*(\a,x)&=&\frac{1}{18900}(\a-1)\a(\a+1)x\l[(2\a-1)(2\a+1)x^4-10(2\a-1)x^2+15\r]\nn
c_8^*(\a,x)&=&\frac{1}{907200}(\a-1)\a(\a+1)\l[(2\a-3)(2\a-1)(2\a+1)x^6+{}\r.\nn
& &\hspace{2cm}{}\l.{}-15(2\a-3)(2\a-1)x^4
+45(2\a-3)x^2-15\r]\nn
c_9^*(\a,x)&=&\frac{1}{28576800}(\a-2)(\a-1)\a(\a+1)x\l[(2\a-3)(2\a-1)(2\a+1)x^6+{}\r.\nn
& &\hspace{2cm}{}\l.{}-21(2\a-3)(2\a-1)x^4+105(2\a-3)x^2-105\r].\n\eea
So we might guess that
\be\la{cstar}\l\{\ba{l}\ds c_{2i}^*(\a,x)=\frac{4(-1)^{i+1}}{(2i)!}{\a+1 \ch i-1}
\hyp{2}{1}{-i+1,\a+\frac{5}{2}-i}{\frac{1}{2}}{x^2},\;i=1,2,3,\ldots\\ \\
\ds c_1^*(\a,x)=0\\ \\
\ds c_{2i+1}^*(\a,x)=\frac{8(-1)^{i+1}}{(2i+1)!}{\a \ch i-1}(\a+1)x\;
\hyp{2}{1}{-i+1,\a+\frac{5}{2}-i}{\frac{3}{2}}{x^2},\\
\hfill i=1,2,3,\ldots.\ea\r.\ee

Hence, we have found the following conjecture.

\vspace{5mm}

{\bf Conjecture.} {\sl The polynomials $\set{\SGP}$ satisfy a linear differential
equation of the form
$$M\sum_{i=0}^{\infty}c_i(x)y^{(i)}(x)+
(1-x^2)y''(x)-2(\a+1)xy'(x)+n(n+2\a+1)y(x)=0,$$
where the coefficients $\l\{c_i(x)\r\}_{i=0}^{\infty}$ are given by
(\ref{cnul}), (\ref{c}) and (\ref{cstar}).}

\vspace{5mm}

We remark that this conjecture would imply that for nonnegative integer
values of $\a$ the polynomials $\set{\SGP}$ satisfy a linear differential
equation of order $2\a+4$.

\vspace{5mm}

Finally, we note that since
$$\hyp{2}{1}{-i+1,\a+\frac{5}{2}-i}{\frac{1}{2}}{1}=\frac{(-\a-2+i)_{i-1}}
{(\frac{1}{2})_{i-1}},\;i=1,2,3,\ldots$$
and
$$\hyp{2}{1}{-i+1,\a+\frac{5}{2}-i}{\frac{3}{2}}{1}=\frac{(-\a-1+i)_{i-1}}
{(\frac{3}{2})_{i-1}},\;i=1,2,3,\ldots$$
we have
\bea \sum_{i=1}^{\infty}c_{2i}^*(\a,-1)&=&
\sum_{i=1}^{\infty}\frac{4}{(\frac{1}{2})_i(1)_i4^i}\frac{(-\a-1)_{i-1}}
{(i-1)!}\frac{(-\a-2+i)_{i-1}}{(\frac{1}{2})_{i-1}}\nn
&=&2\sum_{i=0}^{\infty}\frac{(-\a-1)_i(-\a-1+i)_i}{i!(\frac{3}{2})_i(2)_i
(\frac{1}{2})_i}(\frac{1}{4})^i
=2\sum_{i=0}^{\infty}\frac{(-\a-1)_{2i}}{(2i)!(3)_{2i}}2^{2i}\n\eea
and
\bea \sum_{i=0}^{\infty}c_{2i+1}^*(\a,-1)&=&
-\sum_{i=1}^{\infty}\frac{8(\a+1)}{(1)_i(\frac{3}{2})_i4^i}\frac{(-\a)_{i-1}}
{(i-1)!}\frac{(-\a-1+i)_{i-1}}{(\frac{3}{2})_{i-1}}\nn
&=&-\frac{4}{3}(\a+1)\sum_{i=0}^{\infty}\frac{(-\a)_i(-\a+i)_i}{i!(2)_i
(\frac{5}{2})_i(\frac{3}{2})_i}(\frac{1}{4})^i\nn
&=&-\frac{4}{3}(\a+1)\sum_{i=0}^{\infty}\frac{(-\a)_{2i}}{(2i+1)!(4)_{2i}}2^{2i}
=2\sum_{i=0}^{\infty}\frac{(-\a-1)_{2i+1}}{(2i+1)!(3)_{2i+1}}2^{2i+1}.\n\eea
In the same way we find
$$\sum_{i=1}^{\infty}c_{2i}^*(\a,1)=
2\sum_{i=0}^{\infty}\frac{(-\a-1)_{2i}}{(2i)!(3)_{2i}}2^{2i}$$
and
$$\sum_{i=0}^{\infty}c_{2i+1}^*(\a,1)=
-2\sum_{i=0}^{\infty}\frac{(-\a-1)_{2i+1}}{(2i+1)!(3)_{2i+1}}2^{2i+1}.$$
This implies that
\bea \sum_{i=1}^{\infty}c_i^*(\a,-1)&=&\sum_{i=1}^{\infty}c_{2i}^*(\a,-1)+
\sum_{i=0}^{\infty}c_{2i+1}^*(\a,-1)\nn
&=&2\sum_{i=0}^{\infty}\frac{(-\a-1)_i}{i!(3)_i}2^i=2\;\hyp{1}{1}{-\a-1}{3}{2}\n\eea
and
\bea \sum_{i=1}^{\infty}c_i^*(\a,1)&=&\sum_{i=1}^{\infty}c_{2i}^*(\a,1)+
\sum_{i=0}^{\infty}c_{2i+1}^*(\a,1)\nn
&=&2\sum_{i=0}^{\infty}\frac{(-\a-1)_i}{i!(3)_i}(-2)^i=2\;\hyp{1}{1}{-\a-1}{3}{-2}.\n\eea

\vspace{5mm}

{\bf Remark.} The proof of the conjecture has been found and will be given in
a forthcoming paper \cite{SGP}.

\vspace{5mm}

{\bf Acknowledgement.} The author wishes to thank Professors Desmond Evans,
Norrie Everitt and Lance Littlejohn for their invitation to come to Cardiff,
where the formulas of the conjecture were found.

Further he wishes to thank his uncle Jan Koekoek for his valuable remarks
after reading the first version of this paper.

\end{document}